\documentclass[10pt, a4paper]{article}

\usepackage{amssymb}
\usepackage{amsthm}
\usepackage{amsmath}
 
\allowdisplaybreaks
    
\newcommand{\R}{\mathbf{R}}

\newcommand{\me}{\mathrm{e}} 
\newcommand{\mi}{\mathrm{i}}
\newcommand{\dif}{\mathrm{d}} 
\newcommand{\Dif}{\mathrm{D}}
\newcommand{\elbow}{\, \rfloor \, }
\newcommand{\mycirc}{\mbox{\footnotesize $\circ$}}
\renewcommand{\Re}{\mathbf{Re}}
\renewcommand{\Im}{\mathbf{Im}}

\newtheorem*{mainthm}{Main Theorem}
\newtheorem{thm}{Theorem}
\newtheorem{lemma}[thm]{Lemma}

\newtheorem{prop}[thm]{Proposition}
\newtheorem*{nonumthm}{Theorem}
\newtheorem*{maincor}{Corollary}

\theoremstyle{definition}
\newtheorem{defn}[thm]{Definition}

\begin{document}
    
\title{Deformations of Minimal Lagrangian Submanifolds \\ with
Boundary}

\author{Adrian Butscher \\ Max Planck Institute for Gravitational
Physics \\ email: \ttfamily butscher@aei-potsdam.mpg.de}

\maketitle

\begin{abstract}
    Let $L$ be a special Lagrangian submanifold of a compact
    Calabi-Yau manifold $M$ with boundary lying on the symplectic,
    codimension 2 submanifold $W$.  It is shown how deformations of
    $L$ which keep the boundary of $L$ confined to $W$ can be
    described by an elliptic boundary value problem, and two results
    about minimal Lagrangian submanifolds with boundary are derived
    using this fact.  The first is that the space of minimal
    Lagrangian submanifolds near $L$ with boundary on $W$ is found to
    be finite dimensional and is parametrized over the space of
    harmonic 1-forms of $L$ satisfying Neumann boundary conditions. 
    The second is that if $W'$ is a symplectic, codimension 2
    submanifold sufficiently near $W$, then under suitable conditions,
    there exists a minimal Lagrangian submanifold $L'$ near $L$ with
    boundary on $W'$.
\end{abstract}

\renewcommand{\thefootnote}{\fnsymbol{footnote}} \footnotetext{
Typeset with LaTeX version 3.14159.  AMS Subject Classification 58J05}


\section{Introduction and Statement of Results}

A minimal Lagrangian submanifold of a symplectic manifold $M$ is at
once minimal with respect to the metric of $M$ and Lagrangian with
respect to the symplectic structure of $M$.  Furthermore, when $M$ is
a Calabi-Yau manifold, Harvey and Lawson showed in their seminal paper
\cite[Section III]{hl1} that minimal Lagrangian submanifolds are also
\emph{calibrated}.  A consequence of this property is that minimal
Lagrangian submanifolds satisfy a relatively simple geometric PDE
(simple relative to the equations of vanishing mean curvature and
symplectic form, which they would satisfy by virtue of minimality and
being Lagrangian separately).  Together, the minimal \emph{and}
Lagrangian conditions lead to the following system of equations: $L
\subset M$ is minimal Lagrangian if and only if
\begin{equation}
    \begin{gathered}
	\left.\Im \left( \me^{\mi \theta} \alpha \right) \right|_{L} =
	0 \\
	\left. \omega \right|_{L} = 0 \, ,
    \end{gathered}
    \label{eqn:minlageqn}
\end{equation}
for some real number $\theta$.  Here, $\omega$ is the symplectic form
of $M$ and $\alpha$ is the canonical, non-vanishing, holomorphic
$(n,0)$-form guaranteed by the Calabi-Yau structure of $M$. 

The \emph{calibration form} defined on $M$ is in this case $\Re \left(
\me^{\mi \theta} \alpha \right)$ and thus $\Re \left(\me^{\mi \theta}
\alpha \right) \big|_{L} = \mathrm{Vol}_{L}$.  The submanifold $L$ is
also referred to as \emph{special Lagrangian with phase angle}
$\theta$ in the literature, while if $L$ is minimal Lagrangian with
phase angle $\theta = 0$ then $L$ is simply called \emph{special
Lagrangian}.

Harvey and Lawson and others, for example, have exploited the
geometric structure implicit in the calibration condition in order to
tackle questions related to the existence of minimal Lagrangian
submanifolds.  Harvey and Lawson themselves produce several examples
of minimal Lagrangian submanifolds and give certain general
constructions of such objects.  More recently, Schoen and Wolfson
\cite{sandw} have presented another construction based on variational
methods and are investigating the singularities that can arise there,
while Haskins \cite{haskins} has constructed new examples of special
Lagrangian submanifolds and cones.  The topic of singularities of
special Lagrangian has recently seen many further advances, for
examples in papers by Joyce
\cite{joyce4,joyce5,joyce1,joyce7,joyce6,joyce2,joyce3}, stimulated by
recent developments in string theory and mirror symmetry
\cite{morrison}.

Another approach for producing minimal Lagrangian submanifolds
involves studying the \emph{deformations} of a given minimal
Lagrangian candidate $L$ and selecting those deformations of $L$ which
preserve the minimal Lagrangian condition.  A deformation of a
submanifold $L \subset M$ is a $C^{1}$ family of embeddings $f_{t} : L
\rightarrow M$ of $L$, where $f_{0}$ is the canonical embedding.  The
goal of this analysis is to characterize of the space of submanifolds
\emph{near} $L$ which are still minimal and Lagrangian, and is carried
out by analysing the equations satisfied by minimal Lagrangian
submanifolds using perturbative techniques in the form of the 
Implicit Function Theorem.

The first results in this area were obtained by McLean \cite[Section
3]{mclean} and extended by Hitchin \cite{hitchin}.  Using equations
\eqref{eqn:minlageqn}, McLean identified the first order deformations
of a special Lagrangian submanifold in a Calabi-Yau manifold and
developed a method for integrating them.  He used this to characterize
the space of special Lagrangian submanifolds near $L$ according to the
following theorem.

\begin{nonumthm}[McLean, 1996] 
    Let $M$ be a compact, Calabi-Yau manifold.  The space of special
    Lagrangian submanifolds sufficiently near a given candidate $L
    \subset M$ is finite dimensional and is parametrized over the set
    $\mathcal{H} ^{1}(L)$ of harmonic one-forms of $L$.
\end{nonumthm}

McLean's work has been extended to the case where $M$ is symplectic by
Salur \cite{salur}.  The work presented in this paper extends McLean's
result in another direction --- this time to the case of a minimal
Lagrangian submanifold of a Calabi-Yau manifold $M$ with non-empty
boundary.  This will be done by first creating a framework for
incorporating boundary conditions into the minimal Lagrangian
differential equations.  A theorem characterizing the space of minimal
Lagrangian submanifolds with boundary near a given candidate can then
be formulated that is analogous to McLean's result for submanifolds
with empty boundary.

More precisely, the boundary conditions will arise through geometric
restrictions on the deformations of the special Lagrangian
submanifolds, and an object that will be called a \emph{scaffold} will
be used for this purpose.

\begin{defn} \label{defn:scaffold}
    Let $L$ be a submanifold of $M$ with boundary $\partial L$ and
    inward unit normal vector field $N \in \Gamma \big( T_{\partial
    L}L \big)$.  A \emph{scaffold} for $L$ is a smooth submanifold $W$
    of $M$ with the following properties:
    \begin{enumerate}
	\item $\partial L \subset W$;
	\item $N \in \Gamma \big( T_{\partial L} W \big)^{\omega}$
	(here, $S^{\omega}$ denotes the symplectic orthogonal
	complement of a subspace $S$ of a symplectic vector space $V$,
	defined as $S^{\omega} \equiv \{ v \in V \, : \, \omega(v,s) =
	0 \: \forall \, s \in S \}$);
	\item The bundle $(TW)^{\omega}$ is trivial.
    \end{enumerate}
\end{defn}

\noindent \scshape Remarks: \upshape Condition (2) is a transversality
condition that ensures that $JN$ is perpendicular to $W$, where $J$ is
the complex structure of $M$.  Since $W$ is symplectic, $N$ can not be
parallel to $W$.  It seems reasonable to expect the Main Theorem to
hold with condition (2) replaced by unconstrained transversality of
$N, JN$ to the tangent space of $W$ along $\partial L$, but this
weaker assumption leads to technical problems later on.  In
particular, the boundary value problem appearing in the analysis of
the linearized operator in Section 3.2 is the Hodge system with
\emph{oblique} boundary conditions rather than the Hodge system with
Neumann boundary conditions.  Since this BVP is more difficult to deal
with and leads to geometrically less natural results, the Author has
avoided it here.  Furthermore, it is possible that the most
geometrically natural type of scaffold is when $W$ is a \emph{complex}
submanifold of $M$ (i.~e.\ the tangent spaces of $W$ are invariant
under $J$) and this automatically satisfies the transversality
condition (2) \cite{schoen2}.  \smallskip

\noindent \scshape Further remarks: \upshape  Condition (3) will be
used in the sequel in order to make certain constructions on $W$
possible; also, the above definition of a scaffold has already been
used in \cite{me1}.  \medskip

The main theorem to be proved in this paper uses a scaffold to
introduce a boundary condition according to the following statement.

\begin{mainthm}[Boundary Deformation Theorem] Let $L$ be a special
Lagrangian submanifold of a compact Calabi-Yau manifold $M$ with
non-empty boundary $\partial L$ and let $W$ be a symplectic,
codimension two scaffold for $L$.  Then the space of minimal
Lagrangian submanifolds sufficiently near $L$ (in a suitable
$C^{1,\beta}$ sense to be defined later on) but with boundary on $W$
is finite dimensional and is parametrized over the harmonic 1-forms of
$L$ satisfying Neumann boundary conditions.
\end{mainthm}

This theorem, the analogue of McLean's theorem for special Lagrangian
submanifolds with boundary, confines $\partial L$ to move only along
the scaffold $W$, and imposes the boundary condition in the following
way.  If $f_{t} : L \longrightarrow M$ is a boundary-confining
deformation, then $f_{t}(\partial L) \subset W$ for all $t$. 
Consequently, the deformations field $V = \left.  \frac{\dif}{\dif t}
f_{t} \right|_{t=0}$ can not be arbitrary: it must be tangent to $W$
at $\partial L$.

\medskip \noindent \scshape Remark: \upshape Another important
difference between the Boundary Deformation Theorem and its
predecessor is that deformations amongst all \emph{minimal} Lagrangian
submanifolds are allowed here and not just amongst special Lagrangian
submanifolds.  This will turn out to be a necessary ingredient of the
proof.  \medskip

At the end of this paper, the Boundary Deformation Theorem will be
used to prove an existence result for minimal Lagrangian submanifolds
with boundary in $M$.  A corollary will be proved which demonstrates
the existence of minimal Lagrangian submanifolds near $L$ with
boundaries on neighbouring scaffolds.

\begin{maincor}[Scaffold Deformation Theorem]
    Let $L$ be a special Lagrangian submanifold of a Calabi-Yau
    manifold $M$ and let $W$ be a symplectic, codimension two scaffold
    for $L$.  Furthermore, suppose that the topology of $L$ is such
    that its first Betti number $b^{1}(L)$ vanishes (and thus $L$ has
    no non-trivial harmonic one forms with Neumann boundary
    conditions).  Then if $W'$ is any symplectic, codimension two
    submanifold of $M$ that is sufficiently near $W$ in the same sense
    as in the Main Theorem, then there is a minimal Lagrangian
    submanifold $L'$, near $L$ and with boundary on $W'$.
\end{maincor}

The remainder of this paper will be organized in the following manner. 
In Section \ref{sect:formulating}, the boundary value problem
describing minimal Lagrangian submanifolds with boundary on a scaffold
is formulated, and in Section \ref{sect:proof}, the proof of the main
theorem is undertaken by solving this boundary value problem.  Since
the Implicit Function Theorem is to be used for this purpose, the
linearized operator corresponding to the BVP must be calculated there
and shown to be surjective with finite dimensional kernel isomorphic
in a suitable sense to the harmonic 1-forms of $L$.  The corollary of
the Main Theorem is then proved in Section \ref{sect:cor} using the
machinery constructed in the preceding sections.

\section{Formulating the Boundary Value Problem}

\label{sect:formulating}

\subsection{Introduction}

For the remainder of this paper, assume that $L$ is a given, fixed
\emph{special} Lagrangian submanifold with boundary that is contained
in an ambient Calabi-Yau manifold $M$, and that $M$ possesses a metric
$g$, a symplectic form $\omega$, and compatible complex structure $J$. 
Denote by $\alpha$ the canonical, holomorophic, non-vanishing
$(n,0)$-form of $M$.  Furthermore, assume that $L$ is connected; the
results for non-connected $L$ follow simply by considering each
component of $L$ separately.  The equations \eqref{eqn:minlageqn}
satisfied by minimal Lagrangian submanifolds suggest the definition of
a map whose zero set corresponds to the minimal Lagrangian
submanifolds near $L$.  Suppose that the dimension of $L$ is $n$.  Let
$\mathit{Emb}(L,M)$ denote the set of embeddings of $L$ into $M$
(worry about regularity later) and denote by $\Omega^{k}(L)$ the
$k$-forms of $L$.  Now define $\Phi : \mathit{Emb}(L,M) \times \R
\rightarrow \Omega^{1}(L) \times \Omega^{n}(L)$ by
\begin{equation}
    \label{eqn:map}
    \Phi(f,\theta) = \big( f^{\ast} \omega, f^{\ast} \Im (\me^{\mi
    \theta} \alpha) \big) \, .
\end{equation}
Since $L$ itself is special Lagrangian, $\Phi(i_{L}, 0) = (0,0)$,
where $i_{L}$ is the canonical embedding of $L$.  Another minimal
Lagrangian embedding of $L$, with calibration angle $\theta$, is an
embedding $f$ satisfying $\Phi(f, \theta) = (0,0)$.

The main theorem of this paper consists of finding those embeddings of
$L$ near $i_{L}$ which satisfy $\Phi(f, \theta) = (0,0)$ for some
$\theta$ by means of the Implicit Function Theorem.  The precise
version of the theorem that will be employed is the following.

\begin{thm}[Implicit Function Theorem]
    \label{thm:ift}
    Let $F : B \rightarrow Z$ be a $C^{1}$ map of Banach spaces with
    $F(0) = 0$.  Suppose that there exist closed Banach subspaces $X$
    and $Y$ of $\mathcal{B}$ so that $\mathcal{B} = X \oplus Y$.  If
    $\Dif_{X} F(0)$ is bijective, then there is a neighbourhood
    $\mathcal{U}$ of $0$ in $Y$ and a $C^{1}$ map $\phi : \mathcal{U}
    \rightarrow X$ so that $\phi(0) = 0$ and $F \big( y + \phi(y)
    \big) = 0$ for all $y \in \mathcal{U}$.
\end{thm}

\noindent See \cite[Section 2.5]{amr} for an excellent discussion of
this theorem as well as its proof.  The Implicit Function Theorem thus
provides families of solutions of the equation $F(b) = 0$,
parametrized over the Banach subspace which complements the subspace
on which the linearization of $F$ at $0$ is bijective.  Note that in
the case where $\Dif F(0)$ is surjective with finite dimensional
kernel $K$, then the Implicit Function Theorem holds with $Y = K$ and
$X$ equal to any Banach subspace, necessarily closed, that complements
$K$.  The main theorem will be proved using this special case, while
the corollary will be proved using the more general statement of the
Implicit Function Theorem.

The map $\Phi$, as defined above, does not yet involve Banach spaces. 
Thus in order to apply the Implicit Function Theorem to $\Phi$, a
sufficiently large class of embeddings of $L$ near $i_{L}$ must be
parametrized over a Banach space, and the equation $\Phi(f,\theta) =
(0,0)$ must be solved in this Banach space.  An added difficulty is
that the elements of the Banach space must satisfy a boundary
condition which ensures that $\Phi$ acting on the Banach space is
elliptic.  

\subsection{Imposing Boundary Conditions with a Scaffold}

\label{sect:imposing}

In order to understand why boundary conditions must be imposed on the 
deformations of $L$, one must consider the linearization of the 
operator $\Phi$ at the point $(i_{L},0)$.  

\begin{prop}
    \label{prop:linop}
    Let $\Phi: \mathit{Emb}(L,M) \times \R \rightarrow \Omega^{1}(L)
    \times \Omega^{n}(L)$ be the operator defined in \eqref{eqn:map}. 
    The linearization of $\Phi$ at the point $(i_{L},0)$ is given by
    \begin{equation}
	\label{eqn:linop}
	\Dif \Phi(i_{L}, 0)(V, a) = \big( \dif \eta, \dif \star \eta +
	a \mathrm{Vol}_{L} \big) \, ,
    \end{equation}
    where $V$ is a vector field defined on $L$, $a$ is a real number
    and $\eta = i_{L}^{\ast} ( V \elbow \omega )$.
\end{prop}

\begin{proof}
    
Let $f_{t}:L \rightarrow M$ be a family of embeddings with $f_{0} =
i_{L}$ and $\frac{\dif}{\dif t} f_{t} \big|_{t=0} = V$; and let
$a_{t}$ be a family of real numbers with $a_{0} = 0$ and
$\frac{\dif}{\dif t} a_{t} \big|_{t=0} = a$.  Now,
\begin{equation*}
    \Dif \Phi(i_{L},0) (V,a) = \left. \frac{\dif}{\dif t} \Phi(f_{t}, 
    a_{t}) \right|_{t=0} \, .
\end{equation*}
The calculation of the derivative of $\Phi$ in the $f_{t}$ direction
has already been carried out by McLean in \cite{mclean}.  It remains
only to perform the calculation in the $a_{t}$ direction.  This can be
done by differentiating
\begin{align*}
    \left.  \frac{\dif}{\dif t} \Phi(0,ta) \right|_{t=0} &= \left(0,
    -\left.  \frac{\dif}{\dif t} \Im \left( \me^{-\mi t a} \alpha
    \right) \right) \right|_{t=0} \\
    &= \big( 0, i^{\ast} \big( \Im \left( \mi a \alpha \right) \big)
    \big) \\
    &= \big( 0, a \, i^{\ast} \big( \Re (\alpha) \big) \big) \\
    &= \big( 0, a \, \mathrm{Vol}_{L} \big) \, ,
\end{align*}
by definition of a calibration form.  This calculation, in combination
with McLean's result, completes the proof of the proposition.
\end{proof}

The reason boundary conditions are necessary is that the Hodge
operator $\eta \mapsto (\dif \eta, \dif \star \eta)$ is \emph{not}
elliptic unless it acts upon a space of differential 1-forms
satisfying certain boundary conditions.  From the Hodge theory on
manifolds with boundary \cite{schwarz}, it is known that one such
boundary condition is the \emph{Neumann boundary condition}: the Hodge
operator is elliptic when acting on forms $\eta$ which satisfy
$\eta(N) = 0$ along $\partial L$, where $N$ is the unit normal vector
field of $\partial L$ in $L$.  In the case under consideration here,
$\eta$ arises as the 1-form associated to a deformation of a special
Lagrangian submanifold, and is thus of the form $\eta = V \elbow
\omega$, where $V = \frac{\dif}{\dif t} f_{t} \big|_{t=0}$ is the
corresponding deformation vector field.  The Neumann boundary
condition thus translates into the condition
\begin{equation}
    i^{\ast}_{L} (V \elbow \omega) (Y) = 0 \qquad \Leftrightarrow
    \qquad \omega(V, Y) = 0
    \label{eqn:geobdcond}
\end{equation}
on $V$.  The following proposition shows that this boundary condition
arises naturally if the deformations of $L$ that are considered force
the boundary of $L$ to remain on a scaffold as described in the
introduction.

\begin{prop} \label{prop:symplscaffold}
    Let $L$ be a special Lagrangian submanifold of $M$ and let $W$
    be a scaffold for $L$.  In addition, suppose $W$ is scaffold for
    $L$ which is also symplectic and has codimension 2 in $M$.  Let
    $f_{t} : L \longrightarrow M$ be any deformation of $L$
    satisfying $f_{t} (\partial L) \in W$ for all $t$.  Then the
    deformation vector field $V = \left.  \frac{\dif}{\dif t} f_{t}
    \right|_{t=0}$ corresponding to $f_{t}$ satisfies the elliptic
    boundary condition $\omega(V, N) = 0$, where $N$ is the unit
    normal vector field of $\partial L$ in $L$.
\end{prop}

\begin{proof}

The vector field $V$ must be parallel to $W$ along $\partial L$ as
indicated in the introduction.  But according to the definition of a
scaffold, $N \in (T_{x}W)^{\omega}$ for every $x \in \partial L$. 
Therefore $\omega(N, V) = 0$.
\end{proof}

\subsection{Constructing Scaffold Preserving Deformations}

\label{sect:constructing}

In the proof of McLean's Theorem, deformations of $L$ are parametrized
over the Banach space of $C^{1,\beta}$ sections of the the normal
bundle of $L$ using the exponential map.  That is, for every section
$V$ of the normal bundle of $L$, the exponential map defines an
embedding of $L$ via $\exp(V) : L \rightarrow M$.  Exponential
deformations are, however, not suitable for the present purpose,
because in general $\exp(V)(\partial L)$ will not lie on $W$ because
$W$ is in general not totally geodesic.  Indeed, if $p \in \partial L$
and the geodesic starting at $p$ and heading in the direction of $V$
does not lie in $W$, then $\exp(V)(p) \not\in W$.  Another means of
deforming $L$ is thus necessary if $\partial L$ is to remain confined
to the scaffold under deformation.  One way to avoid the difficulty
described above is to consider the exponential map of a
\emph{different} metric $\hat{g}$ --- one in which $W$ \emph{is}
totally geodesic.  The normal bundle of $L$ with respect to this new
metric, denoted by $\hat{N}L$, will then be used to parametrize
submanifolds near $L$ with boundary on $W$.

Before the metric $\hat{g}$ can be constructed, a lemma concerning the
local structure of $W$ near $\partial L$ is needed.  This is
essentially a version of the Lagrangian Neighbourhood Theorem
\cite[page 99]{mcduff} that is valid for Lagrangian submanifolds with
boundary.

\begin{lemma} \label{lemma:scaffolddarboux}
    Let $W$ be a symplectic submanifold of codimension 2 in $M$ and
    suppose that $L$ is a Lagrangian submanifold with boundary
    $\partial L \subset W$.  Then there exists a tubular neighbourhood
    $\mathcal{U}$ of the boundary and a symplectomorphism $\psi :
    \mathcal{U} \longrightarrow T^{\ast}(\partial L \times \R)$ with
    the following properties:
    \begin{enumerate}
	\item $\psi \big( W \cap \mathcal{U} \big) \subset T^{\ast}
	(\partial L) \times \{0,0\}$;
	\item $\psi(\partial L) = \partial L \times \{0,0\}$; 
	\item $\psi \big( L \cap \mathcal{U} \big) \subset \partial L
	\times \R_{+} \times \{0\}$; and
	\item Let $E$ be any non-zero section of $(TW)^{\omega}$ and
	denote by $(s^{1}, s^{2})$ the coordinates of the $\R^{2}$
	factor.  Then $\psi$ can be constructed so that $\psi_{\ast}
	(E) = \frac{\partial} {\partial s^{1}}$.
    \end{enumerate}
    Here, $T^{\ast}(\partial L \times \R)$ has been identified with 
    $T^{\ast}(\partial L) \times \R^{2}$.
\end{lemma}

\begin{proof}
    
Because $W$ is symplectic, the symplectic form $\omega \big|_{W}$
makes $W$ a symplectic manifold in its own right.  Since $\partial L$
is an isotropic submanifold of $M$ with respect to $\omega$, it is
a compact Lagrangian submanifold of $W$ with respect to $\omega
\big|_{W}$.  Consequently, the usual Lagrangian Neighbourhood Theorem
can be applied to $\partial L$ as a submanifold of $W$ to produce a
neighbourhood $\mathcal{U}_{0}$ and a symplectomorphism $\psi_{0} :
\mathcal{U}_{0} \longrightarrow T^{\ast} (\partial L)$.  The desired
symplectomorphism $\psi$ will be found by extending $\psi_{0}$ off $W$
in a suitable way.

The Symplectic Neighbourhood Theorem \cite[page 98]{mcduff} will be
used to complete the extension.  The theorem applies to two symplectic
manifolds $(M_{1}, \omega_{1})$ and $(M_{2}, \omega_{2})$ containing
symplectic submanifolds $W_{1}$ and $W_{2}$ respectively.  It states
that if there exists a symplectic vector bundle isomorphism $\Psi : (T
W_{1})^{\omega} \longrightarrow (T W_{2})^{\omega}$ that covers a
symplectomorphism $\psi : W_{1} \longrightarrow W_{2}$, then there
exist neighbourhoods $\mathcal{U}_{1}$ and $\mathcal{U}_{2}$ of
$W_{1}$ and $W_{2}$ respectively, along with a symplectomorphism
$\psi_{e} : \mathcal{U}_{1} \longrightarrow \mathcal{U}_{2}$ that
extends $\psi$ (that is, $\left.  \psi_{e} \right|_{W_{1}} = \psi$).

Let $M_{1} = M$, $W_{1} = W$, $M_{2} = T^{\ast} (\partial L) \times
\R^{2}$ and $W_{2} = T^{\ast} (\partial L) \times \{0,0\}$.  Let
$s^{1}$ and $s^{2}$ be the coordinate functions in the $\R^{2}$
factor.  One of the defining conditions for a scaffold is that its
sympectic normal bundle $(TW)^{\omega}$ is trivial.  Hence it is
possible to choose two vector fields $E$ and $F$ which span
$(TW)^{\omega}$ and satisfy $\omega(E, F) = 1$.  Extend this basis to
the neighbourhood $\mathcal{U}_{0}$ and continue to denote the
extended vector fields by $E$ and $F$.  Define an isomorphism $\Psi:
(TW)^{\omega} \longrightarrow \R^{2}$ of symplectic vector bundles by:
\begin{equation}
    \Psi(E_{x} ) = \frac{\partial}{\partial s^{1}}_{(\psi_{0} (x),
    0,0)} \qquad \mbox{and} \qquad \Psi(F_{x} ) =
    \frac{\partial}{\partial s^{2}}_{(\psi_{0} (x), 0,0)}
    \label{eqn:sympnbthm}
\end{equation}
at any $x \in \mathcal{U}_{0}$.  This clearly covers the
symplectomorphism $\psi_{0}$ and is a symplectic map.

The Symplectic Neighbourhood Theorem can now be invoked to yield a
symplectomorphism $\psi_{1}$ extending $\psi_{0}$ between some tubular
neighbourhood of $\partial L$ and a neighbourhood of $\partial L
\times \{0,0\}$ in $T^{\ast} (\partial L) \times \R^{2}$. 
Consequently, $\psi_{1}(W \cap \mathcal{U}_{0}) \subset T^{\ast}
(\partial L)$ and $\psi_{1}(\partial L) = \partial L \times \{0,0\}$. 
The third requirement on the symplectomorphism has not yet been met,
however, but the desired property can be obtained by composing with a
suitable symplectomorphism that acts as a translation in the
transverse Lagrangian directions to $L$.
\end{proof}
    
The Darboux coordinates adapted to $L$ guaranteed by Lemma
\ref{lemma:scaffolddarboux} can be used to construct the metric
$\hat{g}$.  The vector field $E$ should be chosen to be the unit
normal vector field $N$ of $\partial L$ in this case.  The
construction of $\hat{g}$ is accomplished in three separate steps.

\smallskip \noindent \scshape Step 1: \upshape Let $\mathcal{U}$ be
the tubular neighbourhood of $\partial L$ provided by Lemma
\ref{lemma:scaffolddarboux} and $\psi : \mathcal{U} \longrightarrow
T^{\ast} (\partial L) \times \R^{2}$ the symplectomorphism.  Suppose
$s^{1}$ and $s^{2}$ are the Darboux coordinates for the $\R^{2}$
factor of the direct product, and furthermore, one can suppose that
$\psi_{\ast} \big( \frac{\partial}{\partial s^{1}} \big) = N$.  Now
define the metric $g_{1}$ at the point $(x, y, s^{1}, s^{2}) \in
T^{\ast} (\partial L) \times \R^{2}$ as follows:
\begin{subequations} \label{eqn:hatmetric}
\begin{equation} \label{eqn:hatmetric1}
    g_{1}(x,y,s^{1},s^{2}) = (\psi^{-1})^{\ast} \big( \left.  g
    \right|_{W} (\psi(x,y,0,0)) \big) + \dif s^{1} \otimes \dif
    s^{1} + \dif s^{2} \otimes \dif s^{2} \, .
\end{equation}

\smallskip \noindent \scshape Step 2: \upshape Without loss of
generality, the form \eqref{eqn:hatmetric1} can be taken for an entire
tubular neighbourhood $\mathcal{U}_{1}$ of $W$.  This is because the
topological assumption made on $W$ --- that $W$ has trivial normal
bundle --- is enough to guarantee the extension of the coordinates
$s^{1}$ and $s^{2}$ to the entire tubular neighbourhood.  The crucial
difference between the present coordinates and the ones used in Step 1
is that the new coordinates are not necessarily symplectic
\emph{everywhere} (but they remain symplectic near $\partial L$).

\smallskip \noindent \scshape Step 3: \upshape Let $\eta: M
\longrightarrow \R$ be a positive, $C^{\infty}$ cut-off function which
equals 1 inside a tubular neighbourhood $\mathcal{U}_{1}'$ of
$\partial L$ contained in $\mathcal{U}_{1}$, and equals $0$ outside
$\mathcal{U}_{1}$.  Now define the metric $\hat{g}$ by:
\begin{equation}
    \hat{g} = \eta g_{1} + (1-\eta) g \, .
\end{equation}
\end{subequations}

\smallskip It remains to verify that the metric $\hat{g}$ brings about
the following properties.  First, it must be true that if $V$ is any
section of the $\hat{g}$-normal bundle $\hat{N}L$ satisfying the
boundary condition $\omega(V, N) = 0$, then $V$ must be tangent to
$W$; this encodes the boundary condition.  Second, $W$ must be
$\hat{g}$-totally geodesic (at least in a neighbourhood of $\partial
L$) so that $\exp(V)(\partial L) \subset W$ whenever $V$ is
sufficiently small.  The following two propositions deal with these
issues.

\begin{prop} \label{prop:totalgeo}
    The submanifold $W$ is totally geodesic with respect to the metric
    $\hat{g}$ constructed in equations \eqref{eqn:hatmetric}.
\end{prop}

\begin{proof}

Let $\frac{\partial} {\partial z^{1}}, \ldots, \frac{\partial}
{\partial z^{2n-2}}$ be a set of local coordinate vector fields for
the relatively open neighbourhood $W \cap \mathcal{U}_{1}'$.  Then
$\frac{\partial} {\partial z^{1}}, \ldots, \frac{\partial} {\partial
z^{2n-2}}, \frac{\partial} {\partial s^{1}}, \frac{\partial}{\partial
s^{2}}$ determines a set of local coordinate vector fields for the
neighbourhood $\mathcal{U}_{1}'$.  In these coordinates,
\begin{equation*}
    \hat{g}=
    \begin{pmatrix}
	\big( \left.  g \right|_{W} (z) \big)_{ij} & 0 \\
	0 & \delta_{ij}
    \end{pmatrix} \, .
\end{equation*}
Now, 
\begin{align*}
    \left\langle \nabla_{\frac{\partial} {\partial z^{i}}}
    \frac{\partial} {\partial z^{j}} , \frac{\partial} {\partial 
    s^{k}} \right\rangle &= \frac{1}{2} \left( \hat{g}_{z^{i} s^{k}, 
    z^{j}} + \hat{g}_{z^{j} s^{k}, z^{i}} - \hat{g}_{z^{i} z^{j}, 
    s^{k}} \right) \\
    &= -\frac{1}{2} \frac{\partial} {\partial s^{k}} \left( \left.  g
    \right|_{W} \right)_{ij} \\
    &= 0 \, .
\end{align*}
This implies that the second fundamental form of $W$ with respect 
to $\hat{g}$ vanishes; and this, in turn, is equivalent to the fact 
that $W$ is totally geodesic.
\end{proof}

\begin{prop} \label{prop:hatellbdcond}
    Let $L$ be a special Lagrangian submanifold with boundary on the
    symplectic scaffold $W$ and let $N$ be the unit normal vector
    field of $\partial L$ in $L$.  Construct the metric $\hat{g}$
    according to equations \eqref{eqn:hatmetric}.  Suppose $V$ is a
    section of $\hat{N}L$ that satisfies the boundary condition
    $\omega(V, N)=0$.  Then $V$ is tangent to $W$ over $\partial L$.
\end{prop}

\begin{proof}
    
Choose a point $x$ in $\partial L$ and Darboux coordinates at $x$ as
in the constructions above.  Furthermore, assume that $T_{x}\partial
L$ is spanned by $\frac{\partial} {\partial z^{1}}, \ldots,
\frac{\partial} {\partial z^{n-1}}$ and that $\frac{\partial}
{\partial z^{n}}, \ldots, \frac{\partial}{\partial z^{2n-2}}$ are
orthogonal to these vectors.  Since $N$ equals $\frac{\partial}
{\partial s^{1}}$ in these coordinates, it is now easy to see that the
$\hat{g}$-normal bundle of $L$ at $x$ is spanned by the vectors
\begin{equation*}
    \frac{\partial} {\partial z^{n}}, \ldots, \frac{\partial}{\partial
    z^{2n-2}} \quad \mbox{and} \quad \frac{\partial}{\partial 
    s^{2}} + \lambda \frac{\partial}{\partial s^{1}} \, ,
\end{equation*}
for some $\lambda \in \R$.  So, if $V \in \hat{N}_{x}L$ and $\big(
\sum \dif z^{i} \wedge \dif z^{n-1+i} + \dif s^{1} \wedge \dif s^{2}
\big) (V, N) = 0$, then clearly the $\frac{\partial}{\partial s^{2}} +
\lambda \frac{\partial}{\partial s^{1}}$ component of $V$ must vanish
and as a result, $V \in T_{x}W$.
\end{proof}

From elementary metric geometry, it is known that $\widehat{\exp}$ is
a local diffeomorphism on $\hat{N}L$.  Thus the conclusion to be drawn
from Proposition \ref{prop:totalgeo} and Proposition
\ref{prop:hatellbdcond} is that sufficiently small
$\hat{g}$-exponential deformations of sections of $\hat{N}L$
satisfying the boundary condition imposed by the scaffold $W$ are in
one-to-one correspondence with submanifolds near $L$ with boundary on
$W$ that project onto $L$ via $\hat{g}$-nearest point projection. 
This parametrization will now be used to apply the Implicit Function
Theorem to the problem of selecting those sections $V \in \hat{N}L$
which give rise to minimal Lagrangian submanifolds near $L$ with
boundary on $W$.

\section{Proof of the Main Theorem}

\label{sect:proof}

\subsection{Defining the Differential Operator}

The apparatus created in the previous section for deforming the
special Lagrangian submanifold $L \subset M$ such that its boundary
remains confined to the scaffold $W$ can now be used to set up a
differential equation whose solutions correspond to minimal Lagrangian
submanifolds near $L$ with boundary on $W$.  Construct the metric
$\hat{g}$ and the $\hat{g}$-normal bundle $\hat{N}L$ of $L$ as in the
previous section, let $N$ denote the unit normal vector field of
$\partial L$, and define the Banach space
$$\mathcal{X} = \left\{ V \in C^{1,\beta} \big( \Gamma ( \hat{N}L )
\big) : \omega(V, Y) = 0 \right\}
$$
of $C^{1,\beta}$, $\hat{g}$-normal vector fields satisfying the
Neumann boundary condition imposed by $W$.  The notation used here is
the following.  If $B$ denotes any bundle over $L$, then $\Gamma(B)$
denotes the sections of $B$ and $C^{k,\beta} \big( \Gamma(B) \big)$
denotes the set of sections whose $k$ covariant derivatives exist and
are bounded in the $C^{k,\beta}$ norm, which is given by
$$\vert u \vert_{C^{k,\beta}} = \sum_{i=0}^{k} \Vert \nabla^{i} u 
\Vert_{0} + [ \nabla^{k} u ]_{\beta} \, ,$$
for any section $u \in \Gamma(B)$, where $\Vert u \Vert_{0}$ is the
supremum norm of a section $u$ over $L$ and $[u ]_{\beta}$ is its
H\"older coefficient.

Next, denote by $\dif \Omega^{k}(L)$ the set of \emph{exact}
$(k+1)$-forms of $L$ and define the operator
$$\Phi : \mathcal{X} \times \R \rightarrow C^{0,\beta} \big( \dif \,
\Omega^{1}(L) \big) \times C^{0,\beta} \big( \dif \, \Omega^{n-1}(L)
\big)$$
by
\begin{equation}
    \label{eqn:defop}
    \Phi(V, \theta) = \big( \widehat{\exp}(V) \big)^{\ast} \big(
    \omega, -\Im ( \me^{- \mi \theta} \alpha) \big) \, ,
\end{equation}
where $\widehat{\exp}$ is the $\hat{g}$-exponential map as defined in
the previous section.  Note that elements of $C^{0,\beta} \big( \dif
\, \Omega^{k}(L) \big)$ are necessarily of the form $\dif \eta$ for
some $\eta \in C^{1,\beta} \big( \Omega^{k-1}(L) \big)$ by the
Poincar\'e Lemma and the basic properties of H\"older spaces.

Since $L$ is special Lagrangian and $\widehat{\exp}(0) = i_{L}$ (the
standard embedding of $L$) then $\Phi(0,0) = (0,0)$.  
The proof of the Main Theorem of this paper consists of showing that
the Implicit Function Theorem, stated in Theorem \ref{thm:ift}, can be
applied to the operator $\Phi$ in order to find solutions of the
equation $\Phi(V,\theta) = (0,0)$.

\medskip \noindent \scshape Note: \upshape The range of $\Phi$ is
indeed the set of \emph{exact} 1- and $n$-forms.  This is any element 
of the range is homotopic to the zero 1- and $n$-forms, and exactness 
is preserved under homotopy.

\subsection{Analysis of the Linearized Operator}

In order to apply the Implicit Function Theorem to the map $\Phi$ in
the vicinity of the point $(0,0)$, it is necessary to show that $\Phi$
is a $C^{1}$ map of Banach spaces and the linearization $\Dif
\Phi(0,0)$ is bounded and surjective, and that its kernel is
isomorphic to the finite dimensional set of harmonic 1-forms of $L$
that satisfy Neumann boundary conditions.

The continuous differentiability of $\Phi$ as a Banach space map is
straightforward.  Recall now the expression of the linearization of
the minimal Lagrangian equations from Proposition \ref{prop:linop};
since $\frac{\dif}{\dif t} \widehat{\exp}(tV) \big|_{t=0} = V$,
\begin{equation*}
    \Dif \Phi(0,0)(V,a) = \Big( \dif \eta , \dif \star \eta + a
    \mathrm{Vol}_{L} \Big) \, ,
\end{equation*}
where $\eta = i_{L}^{\ast} (V \elbow \omega)$ as before.  This is
clearly a bounded operator, and due to its relative simplicity,
surjectivity is easy to verify.

\begin{prop}
    The operator $\Dif \Phi(0,0) : \mathcal{X} \times \R \rightarrow
    C^{0,\beta} \big( \dif \, \Omega^{1}(L) \big) \times C^{0,\beta}
    \big( \dif \, \Omega^{n-1}(L) \big)$ is surjective.
    \label{prop:surj}
\end{prop}

\begin{proof}

Let $N$ be the unit normal vector field of $\partial L$ and let
$\alpha \in C^{1,\beta}\big( \Omega^{1}(L) \big)$ and $\beta \in
C^{1,\beta} \big( \Omega^{n-1}(L) \big)$.  Consider the system of
equations $\Dif \Phi(0,0)(V, a) = \big(\dif \alpha, \dif \beta \big)$
in the space $\mathcal{X} \times \R$; or in other words, consider
\begin{equation}
    \begin{gathered}
	\dif \eta = \dif \alpha \\
	\dif  \star  \eta = \dif \beta + a \, \mathrm{Vol}_{L} \\
	\eta(N) = 0 \, .
    \end{gathered}
    \label{eqn:obliquehodge}
\end{equation}

Hodge theory for a manifold $L$ with boundary (see \cite{schwarz} for
a thorough explanation of all the details of this theory) shows that a
$k$-form satisfying the equations
\begin{equation*}
    \begin{gathered}
	\dif \eta = \sigma \\
	\dif   \star   \eta = \tau \\
	\eta(N) = 0
    \end{gathered}
\end{equation*}
and possessing a given degree of H\"older regularity can be found if
and only if the following conditions are met:
\begin{enumerate}
    \item $\dif \sigma = 0$;
    
    \item $ \dif \tau = 0$;
    
    \item $\tau(E_{1}, \ldots, E_{k+1}) \big|_{\partial L} = 0$ for
    any collection of vectors $E_{i}$ tangent to $\partial L$;
    
    \item $\int_{L} \sigma \wedge \star  \lambda = 0$ for every 
    harmonic $(k+1)$-form $\lambda$ of $L$ satisfying Neumann boundary 
    conditions;
    
    \item $\int_{L} \star  \tau \wedge \star  \kappa = 0$ for 
    every harmonic $(k-1)$-form $\kappa$ of $L$ satisfying Neumann 
    boundary conditions.
\end{enumerate}
This list of conditions is given in \cite[page 123]{schwarz}.  Note
that these results are actually only stated for $k$-forms with Sobolev
regularity.  But they extend fairly easily to H\"older regularity by
standard techniques of elliptic theory (as explained in \cite{schoen},
for example).

Because of form of the equations \eqref{eqn:obliquehodge}, only
condition (5) above imposes any restriction on the solvability of
these equations.  Thus a 1-form $\eta$ that solves this system of
equations (and that possesses the correct degree of regularity) can be
found if and only if the integrability condition
$$\int_{L} \dif \beta + a \int_{L} \mathrm{Vol}_{L} = 0$$
can be satisfied.  But since $\int_{L} \mathrm{Vol}_{L} =
\mathit{Vol}(L) \neq 0$, it is possible to choose $a$ equal to
$$a = - \frac{\int_{L} \dif \beta}{\mathit{Vol}(L)} \, .$$
Hence the integrability condition can be met, proving that $(\eta,a)
\mapsto (\dif \eta, \dif \star \eta + a \mathrm{Vol} )$ is surjective. 
This, in turn, implies that $\Dif \Phi(0,0)$ is surjective.
\end{proof}

In order to complete the proof of the Main Theorem, it remains to find
the kernel of the linearized operator.  Suppose that the equations
\begin{equation}
    \begin{gathered}
	\dif \eta = 0 \\
	\dif \star \eta + a \, \mathrm{Vol}_{L} = 0 \\
	\eta(N) = 0 
    \end{gathered}
    \label{eqn:usualhodge}
\end{equation}
are satisfied by a 1-form $\eta$ on $L$ and a real number $a$. 
Integrating the second equation over $L$ yields:
\begin{align*}
    a \mathit{Vol}(L) &= - \int_{L} \dif  \star  \eta \\
    &= - \int_{\partial L} i_{\partial L}^{\ast} (\star  \eta)  \\
    &= - \int_{\partial L} \star \big( \eta(N) \big) \\
    &= 0
\end{align*}
where $i_{\partial L}$ is the standard embedding of $\partial L$ in
$M$.  The calculations above hold by Stokes' Theorem as well as by
the properties of the Hodge star operator at the boundary of $L$
(these properties are derived in \cite[Sections 1.2 and
2.1]{schwarz}).  Hence $a=0$ and $\eta$ satisfies the Hodge system
$\dif \eta = \delta \eta = 0$ with the boundary condition $\eta(N) =
0$.  The solutions of these equations are the harmonic 1-forms with
Neumann boundary conditions.  This is a finite dimensional space of
dimension equal to $b^{1}(L)$.  Again, this result can be found in
\cite[Section 2.6]{schwarz}.

All of the hypotheses required by the Implicit Function Theorem are
thus satisfied by the map $\Phi : \mathcal{X} \times \R \rightarrow
C^{0,\beta} \big( \dif \, \Omega^{1}(L) \big) \times C^{0,\beta} \big(
\dif \, \Omega^{n-1}(L) \big)$.  Thus if 
$$K = \big\{ V \in \mathcal{X} : \Dif \Phi(0,0)(V,0) = (0,0) \big\}$$
is the finite dimensional kernel of $\Dif \Phi(0,0)$, there is a
$C^{1}$ map $f : \mathcal{U} \rightarrow \mathcal{X} \times \R$, where
$\mathcal{U} \subset K$ is a neighbourhood of $0$, that satisfies
$\Phi (f(k)) = 0$ for every $k \in \mathcal{U}$.  This completes the
proof of the Main Theorem.  \hfill \qedsymbol

\section{Deformations of the Scaffold}

\label{sect:cor}

The main theorem answers the question of the existence of minimal
Lagrangian submanifolds with boundary on the scaffold $W$ which are
near the given candidate $L$.  A relatively simple extension of the
theory that has been developed so far can be used to answer the
question of the existence of minimal Lagrangian submanifolds on
\emph{neighbouring} scaffolds.  If $W'$ is a symplectic scaffold near
$W$ and there is a special Lagrangian submanifold $L$ with boundary
$\partial L \subset W$, one asks whether there is a special (or
minimal) Lagrangian submanifold $L'$ near $L$ with boundary $\partial
L' \subset W'$.  There is an affirmative answer to this question and
it is provided by once again by the Implicit Function Theorem.

To prove results about the minimal Lagrangian submanifolds on
neighbouring scaffolds, it is necessary first to parametrize nearby
scaffolds over a Banach space in some way.  The symplectic structure
preserving \emph{Hamiltonian deformations} of $W$ will be used for
this purpose: a procedure will be developed which associates a
time-one Hamiltonian flow to each element of the set of $C^{2,\beta}$
sections of the two-dimensional bundle $(TW)^{\omega}$.

In order to understand the details of this construction, let $X$ be a
$C^{2,\beta}$ section in $\Gamma \big( (TW)^{\omega} \big)$ and
suppose $\mathcal{U}$ is a tubular neighbourhood of $W$ which is
symplectomorphic to $W \times \R^{2}$.  Furthermore, suppose that the
Darboux coordinate vector fields $\frac{\partial}{\partial s^{1}}$ and
$\frac{\partial} {\partial s^{2}}$ of the $\R^{2}$ factor coincide
with the unit normal vector field $N$ and the vector field $JN$,
respectively, over the boundary $\partial L \subset W$.  (The
existence of such coordinates follows from Lemma
\ref{lemma:scaffolddarboux}.)  Write $X$ in these coordinates as:
\begin{equation*}
    X(q) = a^{1}(q) \frac{\partial}{\partial s^{1}} + a^{2}(q)
    \frac{\partial}{\partial s^{2}}
\end{equation*}
where $q \in W$ and the $a^{i}$ are functions of $W$.  Now let $\eta:
M \longrightarrow \R$ be a positive, $C^{\infty}$ cut-off function
equal to zero outside $\mathcal{U}$ and equal to one inside a smaller
tubular neighbourhood of $W$, and define the function $H_{X} : M
\longrightarrow \R$ by:
\begin{equation*}
    H_{X}(q,s) = \eta(q,s) \big( -a^{2}(q) s^{1} + a^{1}(q) s^{2}
    \big)
\end{equation*}
for $(q,s) \in \mathcal{U}$ and make $H_{X}$ equal to zero elsewhere. 
Because the symplectic form of $W \times \R^{2}$ is equal to $\left. 
\omega \right|_{W} + \dif s^{1} \wedge \dif s^{2}$, it is easy to see
that the Hamiltonian vector field associated to $H_{X}$ is equal to
$X$ when $s^{1} = s^{2} = 0$; that is, on the submanifold $W$ itself. 
Finally, let $\phi_{X} : M \longrightarrow M$ denote the
time-one Hamiltonian flow associated to the function $H_{X}$.  By
elementary properties of the flow, it is clear that
\begin{equation}
    \left.  \frac{\dif}{\dif t} \phi_{tX} \right|_{t=0} = J \nabla
    H_{X} \, ,
    \label{eqn:hamisadif}
\end{equation}
and if this quantity is restricted to $W$, then it equals $X$.

The map $X_{q} \mapsto \phi_{X}(q)$ for $X_{q} \in ( T_{q}
W)^{\omega}$ is a local diffeomorphism because equation
\eqref{eqn:hamisadif} implies that its linearization at the zero
section is the identity.  Without loss of generality, one can assume
that it is the tubular neighbourhood $\mathcal{U}$ that is
diffeomorphic to a neighbourhood of the zero section in $( T_{q}
W)^{\omega}$.  Hence, any scaffold $W'$ sufficiently near $W$ and
sufficiently $C^{1}$-regular (to ensure that $W'$ projects onto $W$)
is a Hamiltonian deformation of the form $W' = \phi_{X}(W)$ for some
vector field $X \in \Gamma \big( (TW)^{\omega} \big)$ that is
sufficiently close to the zero section.  The $C^{1,\beta}$ sections of
the bundle $(TW)^{\omega}$ can thus be used to parametrize scaffolds
sufficiently close to $W$.

This parametrization of nearby scaffolds leads to the following
deformation operator.  Define the map $\Phi_{1}: C^{1,\beta} \big(
\Gamma \big( (TW)^{\omega} \big) \big) \times \mathcal{X} \times \R
\longrightarrow C^{0,\beta} \big( \dif \, \Omega^{1}(L) \big) \times
C^{0,\beta} \big( \dif \, \Omega^{n-1}(L) \big)$ by
\begin{equation}
    \Phi_{1}(X, V, \theta) = \big( \phi_{X} \mycirc \,
    \widehat{\exp}(V) \big)^{\ast} \big( \omega, -\Im \big( \me^{-\mi
    \theta} \alpha \big) \big) \, .
    \label{eqn:scaffolddefop}
\end{equation}
If $\Phi_{1}(X,V,\theta) = (0,0)$ then the submanifold $L' = \phi_{X}
\mycirc \, \widehat{\exp}(V) (L)$ is minimal Lagrangian with
calibration form $\Re ( \me^{-\mi \theta} \alpha)$.  Furthermore,
$\partial L'$ is contained in $W' = \phi_{X}(W)$ because the
deformation $\widehat{\exp}(V)$ preserves $W$.  The parametrization
and deformation operator constructed here now lead to the final result
of this paper.

\begin{thm} \label{thm:scafexistence}
    Let $L$ be a special Lagrangian submanifold of a Calabi-Yau
    manifold $M$ whose boundary lies on a symplectic, codimension two
    scaffold $W$.  Furthermore, suppose that the topology of $L$
    forces its first Betti number $b^{1}(L)$ to vanish.  Then if $W'$
    is any symplectic, codimension two submanifold of $M$ that is
    sufficiently near $W$ in the sense that $W'$ can be written as
    $\phi_{X}(W)$ for some $X \in C^{1,\beta} \big( \Gamma \big(
    (TW)^{\omega} \big) \big)$ which is sufficiently small, then there
    is a minimal Lagrangian submanifold $L'$ near $L$ and with
    boundary on $W'$.
\end{thm}

\begin{proof}  

The linearization of $\Phi_{1}$ in the $\mathcal{X} \times \R$
directions remains the operator from equation \eqref{eqn:linop}, and
is thus an isomorphism because the triviality condition $b^{1}(L)=0$
has been assumed.  Therefore, the Implicit Function Theorem implies
that there is an open set $\mathcal{U}$ of $0$ in $C^{1,\beta} \big(
\Gamma \big( (TW)^{\omega} \big) \big)$ and a map $G: \mathcal{U}
\rightarrow \mathcal{X} \times \R$ satisfying
\begin{equation}
    \Phi_{1}(X, G(X)) = (0,0) \, .
    \label{eqn:implfnexist}
\end{equation}
Suppose $G(X) = \big( V(X), \theta(X) \big)$.  Then equation
\eqref{eqn:implfnexist} is equivalent to the statement that the
submanifold
\begin{equation}
    \label{eqn:scafsoln}
    \phi_{X} \mycirc \, \widehat{\exp}(V(X)) (L)
\end{equation}
is minimal Lagrangian, calibrated by the differential form $\Re \big(
\me^{-\mi \theta(X)} \alpha \big)$ and has boundary on the scaffold
$\phi_{X}(W)$ (this is symplectic because $\phi_{X}$ is a
symplectomorphism).

Consequently, if $W'$ is any codimension 2, symplectic submanifold of
the form $\phi_{X}(W)$ with $X \in \mathcal{U}$, then the minimal
submanifold with boundary on $W'$ required to prove the theorem is
simply \eqref{eqn:scafsoln}.
\end{proof}

\bigskip

\noindent \bfseries Acknowledgements: \mdseries I would like to thank
my Ph.D. advisor at Stanford University, Rick Schoen, for his
patience, insight and his confidence in me while I was carrying out
the research for this paper.  I would also like to thank Justin
Corvino and Vin de Silva for their inspirations, ideas, and careful
proofreading.

\newpage

\bibliography{bddef}

\providecommand{\bysame}{\leavevmode\hbox to3em{\hrulefill}\thinspace}
\begin{thebibliography}{10}

\bibitem{amr}
R.~Abraham, J.~E. Marsden, and T.~Ratiu, \emph{Manifolds, {T}ensor {A}nalysis,
  and {A}pplications}, second ed., Springer-Verlag, New York, 1988.

\bibitem{me1}
Adrian Butscher, \emph{Regularizing a singular special {L}agrangian variety},
  Submitted February 2001.

\bibitem{hl1}
Reese Harvey and H.~Blaine Lawson, Jr., \emph{Calibrated geometries}, Acta
  Math. \textbf{148} (1982), 47--157.

\bibitem{haskins}
Mark Haskins, \emph{Constructing special {L}agrangian cones}, math.DG/0005164.

\bibitem{hitchin}
Nigel~J. Hitchin, \emph{The moduli space of special {L}agrangian submanifolds},
  Ann. Scuola Norm. Sup. Pisa Cl. Sci. (4) \textbf{25} (1997), no.~3-4,
  503--515, Dedicated to Ennio De Giorgi.

\bibitem{joyce4}
Dominic Joyce, \emph{Constructing special {L}agrangian $m$-folds in
  $\mathbf{C}^m$ by evolving quadrics}, \\math.DG/0008154. {T}o appear in
  {M}athematische {A}nnalen.

\bibitem{joyce5}
\bysame, \emph{Evolution equations for special {L}agrangian 3-folds in
  $\mathbf{C}^3$}, math.DG/0010036. To appear in the {A}nnals of {G}lobal
  {A}nalysis and {G}eometry.

\bibitem{joyce1}
\bysame, \emph{Lectures on {C}alabi-{Y}au and special {L}agrangian geometry},
  math.DG/0108088.

\bibitem{joyce7}
\bysame, \emph{Ruled special {L}agrangian 3-folds in $\mathbf{C}^3$},
  math.DG/0012060. To appear in the Proceedings of the London Mathematical
  Society.

\bibitem{joyce6}
\bysame, \emph{Singularities of special {L}agrangian fibrations and the
  {S}{Y}{Z}-conjecture}, \\math.DG/0011179.

\bibitem{joyce2}
\bysame, \emph{Special {L}agrangian 3-folds and integrable systems},
  math.DG/0101249.

\bibitem{joyce3}
\bysame, \emph{Special {L}agrangian $m$-folds in $\mathbf{C}^m$ with
  symmetries}, math.DG/0008021.

\bibitem{mcduff}
Dusa McDuff and Dietmar Salamon, \emph{Introduction to {S}ymplectic
  {T}opology}, {S}econd ed., The Clarendon Press Oxford University Press, New
  York, 1998.

\bibitem{mclean}
Robert~C. McLean, \emph{Deformations of calibrated submanifolds}, Comm. Anal.
  Geom. \textbf{6} (1998), no.~4, 705--747.

\bibitem{morrison}
David~R. Morrison, \emph{Mathematical aspects of mirror symmetry}, Complex
  Algebraic Geometry (Park City, UT, 1993), Amer. Math. Soc., Providence, RI,
  1997, pp.~265--327.

\bibitem{salur}
Sema Salur, \emph{Deformations of special {L}agrangian submanifolds}, Commun.
  Contemp. Math. \textbf{2} (2000), no.~3, 365--372.

\bibitem{sandw}
R.~Schoen and J.~Wolfson, \emph{Minimizing volume among {L}agrangian
  submanifolds}, Differential {E}quations: {L}a {P}ietra 1996 (Shatah Giaquinta
  and Varadhan, eds.), Proc. of Symp. in Pure Math., vol.~65, 1999,
  pp.~181--199.

\bibitem{schoen}
Richard Schoen, \emph{Lecture {N}otes in {G}eometric {P}{D}{E}s on
  {M}anifolds}, Course given in the Spring of 1998 at Stanford University.

\bibitem{schoen2}
\bysame, \emph{Private communication}, Will be part of the MSRI Lecture Notes
  for the Clay Mathematics Institute Summer School on the Global Theory of
  Minimal Surfaces of the Summer of 2001.

\bibitem{schwarz}
G{\"u}nter Schwarz, \emph{Hodge {D}ecomposition---{A} {M}ethod for {S}olving
  {B}oundary {V}alue {P}roblems}, Springer-Verlag, Berlin, 1995.

\end{thebibliography}
\bibliographystyle{amsplain}

\end{document}